\documentclass{article}

\usepackage[colorlinks=true,citecolor=blue]{hyperref}
\usepackage{tikz}
\usetikzlibrary{positioning}

  \usepackage[preprint]{neurips_2026}

\usepackage[utf8]{inputenc} 
\usepackage[T1]{fontenc}    
\usepackage{hyperref}       
\usepackage{url}            
\usepackage{booktabs}       
\usepackage{amsfonts}       
\usepackage{nicefrac}       
\usepackage{microtype}      
\usepackage{xcolor}         
\usepackage{amsmath}

\title{Introductory Notes on Learning²}

\author{%
  Sai S.~Siddharth\\
  Department of Mechanical Engineering, Thiagarajar College of Engineering\\
  Department of Aeronautics \& Astronautics, Massachusetts Institute of Technology\\
  \texttt{siddharthaerospace@gmail.com or saisidd@mit.edu} \\
   \And
   Maniarasu Ravi \\
    Department of Mechanical Engineering \\
   Thiagarajar College of Engineering \\
   \texttt{rmumech@tce.edu} \\
}

\begin{document}

\maketitle

\begin{abstract}
Although machine learning can be used to predict the evolution of physical
systems from data, a formulation that learns only the system state at each
time leaves the temporal and dynamical structure of the solution to be
resolved within a broad hypothesis space. We introduce \emph{Learning$^2$},
a representation-level framework that structures this space by coupling a
primary representation to a second representation through a known physical
transformation. The resulting cross-representation constraint restricts the
effective hypothesis space and provides an ante-hoc, physically interpretable
criterion for excluding solutions that satisfy the primary representation
alone. We instantiate Learning$^2$ through \emph{EuLaNet}, an
Eulerian--Lagrangian representation for fluid dynamics. Given the velocity
state $u(\mathbf{x},t)$, EuLaNet constructs its induced Lagrangian flow map
$X(\mathbf{a},t)$ through
$\dot{X}(\mathbf{a},t)=u(X(\mathbf{a},t),t)$, from which material transport
and finite-time deformation are derived. The resulting representation
couples the predicted state to the dynamical consequences it induces,
providing a second consistency criterion beyond state-level agreement.

We formalize this construction through an effective hypothesis space
$\mathcal{H}_{L^2}\subseteq\mathcal{H}$ and define the conditions under
which a consequence representation provides discriminative constraints on
candidate solutions. EuLaNet is implemented as a model-independent
representation module, separating the physical constraint from the
downstream learning architecture. This construction provides an ante-hoc
mechanism for physically interpretable constraint in scientific learning
and offers a basis for developing and evaluating broader classes of
\emph{Learning$^2$} architectures.The implementation is open-sourced, with a research community being built around Learning$^2$ to develop and extend the architecture across scientific domains.
\end{abstract}

\section{Introduction}

Machine learning is increasingly being used to approximate the evolution of
physical systems that are otherwise governed by computationally expensive
partial differential equations. In fluid mechanics, data-driven methods
have been developed for reduced-order modeling, flow prediction, control,
and simulation, with learned models increasingly operating directly on
spatiotemporal flow data \citep{brunton2020machine,raissi2019physics,
karniadakis2021physics}. A central formulation is to learn the evolution of
a discretized physical state,
\begin{equation}
    \hat{u}_{t+\Delta t}=f_{\theta}(u_t),
\end{equation}
where $u_t=u(\mathbf{x},t)$ denotes the flow state on a fixed spatial
domain. This formulation has led to increasingly expressive architectures,
from models designed around discretized differential operators
\citep{long2018pde} to neural operators that learn mappings between function
spaces \citep{lu2021deeponet,li2021fourier,kovachki2023neural}.

For fluid mechanics, however, the field $u(\mathbf{x},t)$ is only one
representation of the dynamics. A flow can also be described by following
material points through the domain. Given a velocity field
$u(\mathbf{x},t)$, the corresponding Lagrangian flow map is
\begin{equation}
    \frac{d\mathbf{X}(\mathbf{a},t)}{dt}
    =
    u(\mathbf{X}(\mathbf{a},t),t),
    \qquad
    \mathbf{X}(\mathbf{a},t_0)=\mathbf{a},
\end{equation}
where $\mathbf{a}$ denotes the initial position of a material point.
Eulerian and Lagrangian descriptions are therefore two representations of
the same flow. The Eulerian description resolves the state over space,
whereas the Lagrangian description resolves how material is transported
through that state. Lagrangian descriptions are particularly natural for
transport, dispersion, mixing, and particle motion
\citep{toschi2009lagrangian,schroder2023lagrangian}.

This distinction has direct consequences for learned fluid models. A model
trained only on Eulerian states learns a mapping between field
representations. A model that also represents material trajectories has
access to a second description whose evolution is determined by the
predicted field itself. The relationship is therefore not one of simply
adding another set of features. The two representations are coupled by the
kinematics of the flow. This provides an architectural opportunity that is
different from increasing network capacity or augmenting the input with
independent measurements.

We explore this opportunity through \emph{EuLaNet}, an
Eulerian--Lagrangian architecture for learning fluid dynamics. EuLaNet
starts from an Eulerian CFD solution and constructs its corresponding
Lagrangian representation through the flow map. The architecture retains
the Eulerian field as the primary spatial representation while introducing
the trajectories generated by that field as a second representation of the
same evolving system. The resulting model can therefore be viewed as
learning a physical state together with a representation of the motion
implied by that state.

The idea is closely related to a broader class of approaches that build
physical structure into learned dynamics. PDE-Net incorporates
differential operators into a neural architecture for learning dynamical
systems \citep{long2018pde}, while physics-informed methods introduce
governing equations and physical constraints into the learning objective
\citep{raissi2019physics,karniadakis2021physics}. Learned simulators have
also explored particle and mesh representations as alternative state spaces
for physical dynamics \citep{sanchezgonzalez2020learning,pfaff2021learning}.
EuLaNet takes a different architectural route by coupling two
representations that are analytically related through the underlying flow.

We refer to this construction as \emph{Learning$^2$}. The notation denotes
learning in two coupled representations of a physical system. In the
present case, the first representation is Eulerian and the second is
Lagrangian. The purpose of Learning$^2$ is not to prescribe a particular
neural architecture or a particular loss function. It is a design principle
for asking whether the learning of one physical representation can be
organized through another representation of the same dynamics. EuLaNet is
our concrete realization of this principle for fluid flow.

This viewpoint is also distinct from recent Eulerian--Lagrangian neural
models. DeepLag combines Eulerian observations with adaptively sampled
Lagrangian particles and uses inferred particle dynamics to improve fluid
prediction \citep{ma2024deeplag}. Our emphasis is instead on the
representation coupling itself. EuLaNet is used to examine what changes
when the Eulerian state and its dynamically induced Lagrangian description
are made components of the same learning architecture. The resulting
questions concern prediction accuracy, the behavior of the learning
process, the stability of the learned dynamics, and the consistency between
the field that is predicted and the trajectories that it generates.

The contribution of this work is therefore architectural. We formalize
Learning$^2$ as a way of constructing scientific learning systems from
multiple dynamically related representations, instantiate it through
EuLaNet for Eulerian--Lagrangian fluid dynamics, and study the resulting
learning behavior against an Eulerian-only formulation. The objective is
not to close the space of possible Learning$^2$ architectures, but to
provide a concrete construction from which that space can be explored.
The Eulerian--Lagrangian pair is one such physically motivated instance.

\section{Learning$^2$}

We introduce Learning$^2$ as an architectural principle for learning
problems in which the target admits multiple mathematically related
representations. The principle is to use one representation not merely as
additional information about the target, but as a structural constraint on
how another representation is learned. This distinguishes Learning$^2$ from
standard representation learning, in which intermediate representations
primarily provide progressively transformed features for the final
prediction \citep{bengio2013representation,lecun2015deep,goodfellow2016deep}.

\subsection{Learning through representation}

Consider a target $y$ inferred from an input $x$. A conventional learning
problem seeks a mapping

\begin{equation}
    \hat{y}=f_{\theta}(x),
\end{equation}

with solutions selected according to a loss defined on the target
representation. Let $\mathcal{R}(y)$ denote another representation of the
same underlying quantity. Learning$^2$ considers architectures in which
the relation between $y$ and $\mathcal{R}(y)$ is incorporated into the
construction of the learning problem.

The distinction is important. An additional representation can be supplied
to a model as another source of information without changing the structure
of the underlying learning problem. Learning$^2$ instead requires the
relationship between representations to constrain the construction of the
solution. The second representation therefore acts on the space of
admissible solutions rather than simply increasing the information
available to the learner.

Let $\mathcal{H}$ denote the hypothesis space associated with the
single-representation problem. The introduction of a coupled
representation defines an effective hypothesis space

\begin{equation}
    \mathcal{H}_{\mathcal{R}}
    =
    \left\{
        f\in\mathcal{H}
        \;\middle|\;
        f \sim \mathcal{R}(y)
    \right\},
\end{equation}

where $\sim$ denotes the mathematical relationship imposed between the
target representation and $\mathcal{R}(y)$. When this relationship excludes
otherwise admissible hypotheses,

\begin{equation}
    \mathcal{H}_{\mathcal{R}}\subseteq\mathcal{H}.
\end{equation}

Learning$^2$ therefore introduces structure into the hypothesis space
through a representation that is coupled to the quantity being learned.
This provides a mechanism by which a learning architecture can reduce its
effective search space without requiring the target itself to be
represented in a lower-dimensional space or the model to have greater
capacity.

\subsection{Architectural consequence}

The effect of the second representation is consequently not restricted to
the final prediction. By constraining the set of admissible solutions, the
representation relation changes the geometry of the learning problem and
can alter the optimization trajectory through that space. This is
consistent with the established role of architectural structure as an
inductive bias, where restrictions imposed by a model architecture determine
which functions can be represented or preferentially learned
\citep{cohen2016inductive,bronstein2021geometric}.

The distinction can be expressed by comparing the two learning spaces

\begin{equation}
    \mathcal{S}_{1}
    =
    \left\{
        \hat{y}
        \;\middle|\;
        \hat{y}\in\mathcal{H}
    \right\},
\end{equation}

and

\begin{equation}
    \mathcal{S}_{2}
    =
    \left\{
        \hat{y}
        \;\middle|\;
        \hat{y}\in\mathcal{H}_{\mathcal{R}}
    \right\}.
\end{equation}

The Learning$^2$ construction is useful when the additional representation
contains structure that is informative about the admissibility of a
solution, rather than merely information that improves prediction. The
resulting gain may therefore arise from the organization of the learning
space itself, including changes in convergence, stability, generalization,
or long-horizon behavior.

This establishes the central distinction of Learning$^2$. The objective is
not to learn multiple representations independently, nor simply to augment
a model with auxiliary features. The objective is to construct a learning
system in which a mathematically related representation participates in
determining the space through which the target representation is learned.

\subsection{Design principle}

Learning$^2$ does not prescribe a particular network, loss function, or
representation. Its defining components are a target representation
$y$, a related representation $\mathcal{R}(y)$, and a mechanism by which
their relationship constrains the learning process. The choice of this
relationship determines the resulting architecture and the structure
introduced into its learning space.

This formulation provides a general design principle for scientific
machine learning. Whenever a physical system admits multiple mathematical
descriptions of the same underlying state, those descriptions need not be
treated as independent modeling choices. One representation can instead be
used to structure the learning of another. In the present work, we test this
principle through EuLaNet by coupling Eulerian and Lagrangian descriptions
of fluid dynamics.

\section{EuLaNet}

EuLaNet instantiates the Learning$^2$ principle for fluid dynamics by
constructing a coupled Eulerian--Lagrangian representation directly from
time-resolved CFD solutions. The representation preserves the complete
Eulerian flow field while introducing a sparse material description derived
from the same velocity field. The two descriptions are connected through
explicit spatial and temporal correspondence rather than being treated as
independent feature sets.

\subsection{Eulerian state}

Let the CFD solution be available at discrete times
${t_n}_{n=0}^{T-1}$ on a possibly time-dependent computational mesh.
At each snapshot, the mesh coordinates and corresponding flow variables are

\begin{equation}
X_n =
\left\{
\mathbf{x}_{i,n}
\right\}_{i=1}^{N},
\qquad
E_n =
\left\{
\mathbf{e}_{i,n}
\right\}_{i=1}^{N},
\end{equation}

where $\mathbf{x}_{i,n}\in\mathbb{R}^{2}$ denotes the physical location of
mesh point $i$ at time $t_n$ and
$\mathbf{e}_{i,n}\in\mathbb{R}^{d_E}$ contains the corresponding CFD state.
The complete Eulerian representation is therefore

\begin{equation}
E\in\mathbb{R}^{T\times N\times d_E},
\qquad
X\in\mathbb{R}^{T\times N\times 2}.
\end{equation}

The representation retains the full CFD state rather than reducing the
solution to a prescribed subset of observables. For the validated NACA0015
case, this gives

\begin{equation}
E\in\mathbb{R}^{150\times44100\times20},
\qquad
X\in\mathbb{R}^{150\times44100\times2}.
\end{equation}

The CFD solutions are generated with SU2, an open-source framework for
PDE analysis and PDE-constrained optimization on unstructured meshes
\citep{economon2016su2}. The use of time-dependent and moving meshes is
consistent with the dynamic-mesh formulations developed within the SU2
framework for unsteady aerodynamic problems
\citep{economon2015unsteady}.

\subsection{Lagrangian material representation}

The second representation is constructed from material probes embedded in
the Eulerian solution. Let

\begin{equation}
P_0=
\left\{
\mathbf{x}_{p,0}
\right\}_{p=1}^{N_p}
\end{equation}

denote the initial probe locations. Each probe is treated as a passive
material point whose trajectory is determined by the Eulerian velocity
field,

\begin{equation}
\frac{d\mathbf{x}_p}{dt}
=
\mathbf{u}(\mathbf{x}_p,t),
\qquad
\mathbf{x}_p(t_0)=\mathbf{x}_{p,0}.
\end{equation}

For the reference configuration, the initial probes are generated from a
$40\times40$ grid and filtered against the computational domain, yielding
$N_p=1184$ valid material probes. The probe set is therefore sparse relative
to the $44{,}100$-point Eulerian field, while remaining embedded in the same
physical domain.

Because the CFD velocity is available only at discrete snapshots, EuLaNet
constructs the velocity used during particle integration by temporal
interpolation. For $t=t_n+\alpha\Delta t$, with
$\alpha\in[0,1]$,

\begin{equation}
\mathbf{u}(\mathbf{x},t)
=
(1-\alpha)\mathbf{u}_n(\mathbf{x})
+
\alpha\mathbf{u}_{n+1}(\mathbf{x}).
\end{equation}

The resulting trajectory is integrated using fourth-order Runge--Kutta.
For a substep $h$,

\begin{align}
\mathbf{k}_1 &=
f(\mathbf{x}_n,t_n),\\
\mathbf{k}_2 &=
f\left(\mathbf{x}_n+\frac{h}{2}\mathbf{k}_1,
t_n+\frac{h}{2}\right),\\
\mathbf{k}_3 &=
f\left(\mathbf{x}_n+\frac{h}{2}\mathbf{k}_2,
t_n+\frac{h}{2}\right),\\
\mathbf{k}_4 &=
f(\mathbf{x}_n+h\mathbf{k}_3,t_n+h),
\end{align}

with

\begin{equation}
\mathbf{x}_{n+1}
=
\mathbf{x}_n+
\frac{h}{6}
\left(
\mathbf{k}_1+2\mathbf{k}_2+2\mathbf{k}_3+\mathbf{k}_4
\right).
\end{equation}

The resulting discrete flow map is

\begin{equation}
\Phi_{t_0}^{t}
:
\mathbf{x}_0
\mapsto
\mathbf{x}(t;\mathbf{x}_0,t_0),
\end{equation}

which records the material position generated by the Eulerian velocity
field. For the reference case, the complete trajectory representation has
shape

\begin{equation}
\Phi\in
\mathbb{R}^{150\times1184\times2}.
\end{equation}

\subsection{Finite-time deformation}

Particle displacement alone describes transport but does not characterize
the deformation of the material neighborhood surrounding a particle. EuLaNet
therefore derives local deformation from the flow map.

For a trajectory initialized at $\mathbf{x}_0$, the deformation gradient is

\begin{equation}
F_{t_0}^{t}(\mathbf{x}_0)
=
\frac{\partial
\Phi_{t_0}^{t}(\mathbf{x}_0)}
{\partial\mathbf{x}_0}.
\end{equation}

For an infinitesimal initial displacement
$\delta\mathbf{x}_0$, the corresponding displacement at time $t$ satisfies

\begin{equation}
\delta\mathbf{x}_t
\approx
F_{t_0}^{t}\delta\mathbf{x}_0.
\end{equation}

The deformation gradient therefore provides a local linearization of the
material transport map. EuLaNet subsequently forms the right
Cauchy--Green deformation tensor

\begin{equation}
C_{t_0}^{t}
=
\left(F_{t_0}^{t}\right)^{T}
F_{t_0}^{t}.
\end{equation}

For an infinitesimal material vector,

\begin{equation}
\left|
\delta\mathbf{x}_t
\right|^2
=
\delta\mathbf{x}_0^{T}
C_{t_0}^{t}
\delta\mathbf{x}_0.
\end{equation}

The eigenvalue problem

\begin{equation}
C_{t_0}^{t}\mathbf{n}_i
=
\lambda_i\mathbf{n}_i
\end{equation}

then provides the principal directions and squared principal stretches of
the material deformation. In two dimensions, the largest eigenvalue
$\lambda_{\max}$ determines the maximum principal stretch

\begin{equation}
\sigma_{\max}
=
\sqrt{\lambda_{\max}}.
\end{equation}

EuLaNet also derives the finite-time Lyapunov exponent over a window
$T=t-t_0$,

\begin{equation}
\operatorname{FTLE}
=
\frac{1}{2|T|}
\log\lambda_{\max}
\left(
C_{t_0}^{t}
\right),
\end{equation}

providing a scalar measure of finite-time material stretching
\citep{shadden2005definition,allshouse2015surfaces}.

For the rolling-window representation used in the dataset, deformation is
computed between the current and future flow-map states. If

\begin{equation}
J_k
=
\frac{\partial\mathbf{x}(t_k)}
{\partial\mathbf{x}_0},
\qquad
J_{k+W}
=
\frac{\partial\mathbf{x}(t_{k+W})}
{\partial\mathbf{x}_0},
\end{equation}

then the deformation over the window is

\begin{equation}
F_{k\rightarrow k+W}
=
J_{k+W}J_k^{-1}.
\end{equation}

This formulation measures deformation relative to the material state at the
beginning of the selected window rather than relative to the original
release configuration.

\subsection{Eulerian--Lagrangian correspondence}

The two representations cannot be combined by simple concatenation.
A Lagrangian particle changes spatial position with time, while the Eulerian
field is indexed by the spatial coordinates of the CFD mesh at each
snapshot. EuLaNet therefore introduces an explicit correspondence operator.

For particle $p$ at snapshot $t_n$, let

\begin{equation}
i^*
=
\arg\min_i
\left|
\mathbf{x}_p(t_n)-\mathbf{x}_{i,n}
\right|_2.
\end{equation}

The corresponding Eulerian state is then associated with the same-time
mesh location

\begin{equation}
\mathbf{x}_p(t_n)
\longleftrightarrow
\mathbf{x}_{i^*,n},
\end{equation}

with correspondence distance

\begin{equation}
d_p(t_n)
=
\left|
\mathbf{x}_p(t_n)-\mathbf{x}_{i^*,n}
\right|_2.
\end{equation}

The correspondence is constructed independently at every CFD snapshot.
Consequently, the mapping is

\begin{equation}
\mathcal{C}:(p,n)
\mapsto
(i^*,n,d_p),
\end{equation}

rather than a fixed spatial mapping from the initial particle position to
the mesh. This distinction is essential for the moving-mesh configuration,
because the physical location of a mesh point changes with time.

Learning² assumes that a meaningful secondary representation and a known transformation between representations are available. Its effectiveness is consequently dependent on the choice of representations and the quality of the underlying physical description. The present work demonstrates the framework through Eulerian–Lagrangian fluid dynamics, and its applicability to other physical systems remains to be established.All reported EuLaNet computations were CPU-only. The implementation was executed on an Intel i5-10300H CPU (4C/8T, 2.50 GHz) with 15.8 GB RAM. Lagrangian flow-map generation for the \(40\times40\) probe grid required approximately 3.4 hours using single-core Python, with peak memory usage below 2 GB. Subsequent EuLaNet dataset construction required less than 1 minute.

\subsection{Complete representation}

The resulting EuLaNet representation consists of three coupled components,

\begin{equation}
\mathcal{D}
=
\left(
E,X,L,\mathcal{C}
\right),
\end{equation}

where $E$ is the complete Eulerian field, $X$ is the time-dependent mesh
geometry, $L$ contains the sparse Lagrangian transport and deformation
state, and $\mathcal{C}$ defines their spatial and temporal correspondence.

For the validated reference case,

\begin{equation}
L\in\mathbb{R}^{148000\times16}.
\end{equation}

The Lagrangian feature vector contains the initial and current particle
positions, displacement, deformation quantities, principal stretches,
stretch directions, FTLE, and validity information. The complete
representation consequently retains both the instantaneous Eulerian state
and finite-time material response rather than replacing one with the other.

This construction is intentionally model-independent. EuLaNet does not
specify whether the resulting representation is consumed by a CNN, GNN,
Transformer, recurrent model, or another learning architecture. Its role is
to expose the physical relation between the two representations in a form
that can subsequently be incorporated into a learning system.

This provides the concrete realization of Learning$^2$ studied in this
work. The Eulerian state defines the field representation, while the
Lagrangian flow map and its derived deformation quantities provide a
mathematically related representation that constrains and organizes the
learning problem through explicit correspondence. The subsequent
experiments test whether this additional structure produces the predicted
reduction in effective search space and corresponding gains in learning
behavior relative to field-only learning.

\section{Learning$^2$ as a Representation Constraint}

EuLaNet provides a concrete instance of a more general architectural
construction. The defining property is not the use of multiple
representations, but the existence of a structured transformation between
them that imposes an additional consistency condition on candidate
solutions.

Let $R_1(U)$ denote the primary representation of a physical state $U$.
A Learning$^2$ construction introduces a second representation through a
transformation
\begin{equation}
    R_2(U)=\mathcal{T}[R_1(U)],
\end{equation}
where $\mathcal{T}$ is a physically or mathematically defined operator.
The transformation must preserve a relationship between the two
representations that is intrinsic to the underlying system. Coordinate
changes, tensor reshaping, arbitrary embeddings, or direct concatenation
do not satisfy this criterion.

The second representation must additionally provide a discriminative
constraint on candidate solutions. Let $\hat{U}_1$ and $\hat{U}_2$ be two
candidate solutions that are approximately equivalent under the primary
representation,
\begin{equation}
    d_1(R_1(\hat{U}_1),R_1(U^\star))
    \approx
    d_1(R_1(\hat{U}_2),R_1(U^\star)).
\end{equation}
A valid Learning$^2$ representation provides a second metric under which
these candidates can be distinguished,
\begin{equation}
    d_2(\mathcal{T}[R_1(\hat{U}_1)],R_2(U^\star))
    \neq
    d_2(\mathcal{T}[R_1(\hat{U}_2)],R_2(U^\star)).
\end{equation}
The second representation therefore acts as a constraint on the effective
solution space rather than as an additional feature channel.

\begin{figure}[t]
    \centering
    \begin{tikzpicture}[
        node distance=7mm and 10mm,
        box/.style={
            draw,
            rounded corners=2pt,
            align=center,
            minimum height=9mm,
            minimum width=31mm,
            font=\small
        },
        arrow/.style={->, semithick}
    ]

        \node[box] (u) {Physical state $U$};

        \node[box, below left=of u] (r1)
            {Primary representation\\$R_1(U)$};

        \node[box, below right=of u] (r2)
            {Consequence representation\\$R_2(U)$};

        \node[box, below=10mm of u] (T)
            {$R_2(U)=\mathcal{T}[R_1(U)]$};

        \node[box, below=of T] (constraint)
            {Cross-representation\\consistency constraint};

        \draw[arrow] (u) -- (r1);
        \draw[arrow] (u) -- (r2);
        \draw[arrow] (r1) -- (T);
        \draw[arrow] (T) -- (r2);
        \draw[arrow] (T) -- (constraint);

    \end{tikzpicture}
    \caption{Learning$^2$ representation constraint. A second
    representation is generated from the primary representation through a
    physically or mathematically defined transformation and provides an
    independent consistency condition on candidate solutions.}
    \label{fig:learning2_constraint}
\end{figure}
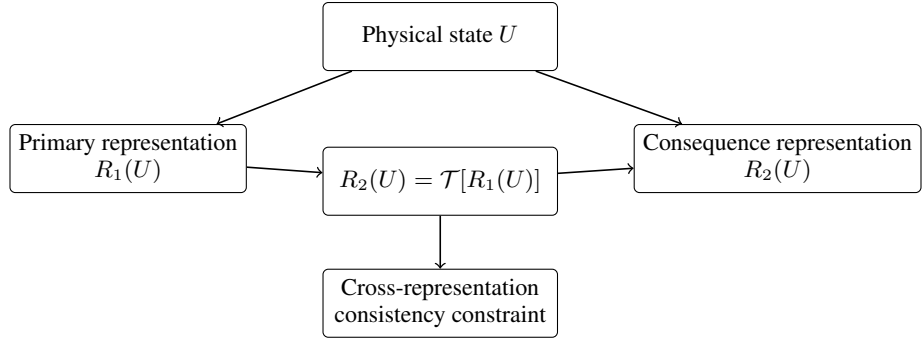

For EuLaNet, the primary representation is the Eulerian velocity field,
\begin{equation}
    R_E(U)=u(\mathbf{x},t),
\end{equation}
and the consequence representation is the Lagrangian flow map,
\begin{equation}
    R_L(U)=X(\mathbf{a},t),
    \qquad
    \frac{dX}{dt}=u(X,t).
\end{equation}
The transformation $\mathcal{T}$ is therefore the material evolution
operator. Finite-time deformation provides further consequences of this
mapping,
\begin{equation}
    u
    \xrightarrow{\mathcal{T}}
    X
    \xrightarrow{\nabla_{\mathbf a}}
    F
    \xrightarrow{F^{\mathsf T}F}
    C
    \xrightarrow{\lambda_{\max}}
    \mathrm{FTLE}.
\end{equation}
These quantities depend on the temporal evolution induced by the field and
cannot be reduced to an independent observation at a single Eulerian
location.

The resulting admissible solution set can be written as
\begin{equation}
    \mathcal{S}_{E}(\epsilon)
    =
    \left\{
        \hat{U}:
        d_E(R_E(\hat{U}),R_E(U^\star))
        \leq \epsilon
    \right\},
\end{equation}
for an Eulerian-only learner, and
\begin{equation}
    \mathcal{S}_{L^2}(\epsilon,\delta)
    =
    \left\{
        \hat{U}\in\mathcal{S}_{E}(\epsilon):
        d_L(\mathcal{T}[R_E(\hat{U})],R_L(U^\star))
        \leq\delta
    \right\}.
\end{equation}
Hence,
\begin{equation}
    \mathcal{S}_{L^2}(\epsilon,\delta)
    \subseteq
    \mathcal{S}_{E}(\epsilon).
\end{equation}
The inclusion describes a restriction of the effective hypothesis space,
not a reduction in neural parameter dimensionality. Candidate fields that
satisfy the primary representation but produce inconsistent dynamical
consequences are excluded by the coupled representation.

This formulation also defines the intended role of EuLaNet within a
larger model architecture. EuLaNet supplies the representation constraint
through Eulerian--Lagrangian correspondence, material transport, and
finite-time deformation. The downstream learner remains independent of
this construction and may be implemented using neural operators, graph
networks, convolutional architectures, sequence models, or other learned
dynamical systems. The representation layer can therefore be treated as a
module rather than as a complete neural architecture.

The resulting design space is defined by the choice of $R_1$, the
transformation $\mathcal{T}$, and the resulting consequence representation
$R_2$. A new Learning$^2$ architecture requires a representation pair for
which the transformation is structurally defined and the second
representation provides discriminative constraints that are not captured
by the primary representation alone. EuLaNet establishes this construction
for Eulerian and Lagrangian fluid dynamics. Extending the framework to
other physical systems reduces to identifying analogous representation
pairs and their corresponding transformations.

Learning$^2$ therefore defines a representation-level constraint on
scientific learning systems. EuLaNet is one realization in which material
dynamics provide that constraint. The resulting formulation provides a
direct basis for constructing new architectures whose objective is not
only to reproduce a physical state, but also to reject solutions whose
induced dynamics are inconsistent with that state. The implementation and accompanying tools are publicly available at the EuLaNet GitHub repository \citep{eulanet2026}, providing a starting point for developing and evaluating new Learning$^2$ architectures.

\small

\bibliographystyle{plainnat}
\bibliography{references}

\appendix

\section{Technical appendices and supplementary material}
\subsection{Experimental Setup}
We provide the complete experimental configuration, including the dataset
construction, train and test cases, temporal sampling, spatial resolution,
and all parameters required to reproduce the reported experiments. The
experiments use the Eulerian--Lagrangian representation described in the
main paper. The complete implementation and data-generation procedure are
available in the EuLaNet repository \citep{eulanet2026}.

\subsection{Computational Resources}
All reported EuLaNet computations were CPU-only. The implementation was
executed on an Intel i5-10300H CPU (4C/8T, 2.50 GHz) with 15.8 GB RAM.
Lagrangian flow-map generation for the $40\times40$ probe grid required
approximately 3.2--3.5 hours using single-core Python, with peak memory
usage below 2 GB. Subsequent EuLaNet dataset construction required less
than 1 minute.

\subsection{Reproducibility}
The EuLaNet implementation, configuration files, dataset-generation
scripts, and instructions required to reproduce the reported representation
are provided in the accompanying repository. The exact command used to
construct the reported $40\times40$ representation is documented there.

\subsection{Limitations and Assumptions}
Learning$^2$ assumes that a meaningful secondary representation and a known
transformation between representations are available. Its effectiveness is
therefore dependent on the choice of representations and the physical
description used to construct the transformation. The present work
demonstrates the framework through Eulerian--Lagrangian fluid dynamics, and
its applicability to other physical systems remains to be established.

\subsection{Statistical Considerations}
Where applicable, experimental variability and repeated runs are reported
with the corresponding results. Since the present study focuses primarily
on the representation and its construction rather than statistical
estimation across a population of independent samples, statistical
significance tests are not used unless otherwise stated.

\subsection{Existing Assets and Licenses}
The implementation uses the documented software dependencies and publicly
available scientific computing tools listed in the repository. Their
respective licenses and versions are retained in the project documentation.
Any external datasets, software packages, or other assets used in the
experiments are identified with their corresponding sources and licenses.

\subsection{New Research Assets}
EuLaNet introduces a new Eulerian--Lagrangian representation and associated
dataset-generation tools. Their construction, assumptions, processing
procedure, and usage are documented in the repository accompanying this
submission.

\subsection{Ethics and Responsible Research}
This work does not involve human participants, human-subject research,
crowdsourcing, personal data, or personally identifiable information. The
experiments use computational fluid dynamics data and derived
Eulerian--Lagrangian representations. No sensitive or human-related data are
collected or processed. The released implementation is intended for
scientific research and is accompanied by documentation of its computational
requirements and usage. To the best of our knowledge, this work does not
introduce foreseeable risks related to privacy, discrimination, surveillance,
deception, or harm to individuals.

\newpage

\end{document}